\documentclass[12pt]{interact}
\usepackage{epstopdf}
\usepackage[caption=false]{subfig}
\usepackage[utf8]{inputenc}
\usepackage{hyperref}

\usepackage{graphicx}
\usepackage{tabularx}
\usepackage{verbatim}
\usepackage{amsmath}
\usepackage[ruled]{algorithm2e}
\SetKwComment{Comment}{$\triangleright$}{}
\usepackage{float}
\usepackage{svg}
\usepackage{tikz}
\usepackage{diagbox}
\usepackage{setspace}
\theoremstyle{plain}

\theoremstyle{definition}

\theoremstyle{remark}

\begin{document}

\title{Identification and prediction of time-varying parameters of COVID-19 model: a data-driven deep learning approach}

\author{
\name{Jie Long\textsuperscript{a}\thanks{Jie Long. Email: jl7u@mtmail.mtsu.edu}, A. Q. M. Khaliq\textsuperscript{b}\thanks{A. Q. M. Khaliq. Email: abdul.khaliq@mtsu.edu} \& K. M. Furati\textsuperscript{c}\thanks{K. M. Furati. Email: kmfurati@kfupm.edu.sa}}
\affil{\textsuperscript{a,b}Department of Mathematical Sciences, Middle Tennessee State University, Murfreesboro, TN 37132, USA; \textsuperscript{b}Department of Mathematics and Statistics, King Fahd University of Petroleum and Minerals, Dhahran 31261, Saudi Arabia}
}

\maketitle

\begin{abstract}
Data-driven deep learning provides efficient algorithms for parameter identification of epidemiology models. Unlike the constant parameters, the complexity of identifying time-varying parameters is largely increased. In this paper, a variant of physics-informed neural network (PINN) is adopted to identify the time-varying parameters of the Susceptible-Infectious-Recovered-Deceased model for the spread of COVID-19 by fitting daily reported cases. The learned parameters are verified by utilizing an ordinary differential equation solver to compute the corresponding solutions of this compartmental model. The effective reproduction number based on these parameters is calculated. Long Short-Term Memory (LSTM) neural network is employed to predict the future weekly time-varying parameters. The numerical simulations demonstrate that PINN combined with LSTM yields accurate and effective results.
\end{abstract}


\section{Introduction}

The novel SARS-CoV-2 (severe acute respiratory syndrome coronavirus~2) has spread all over the world since its discovery at the end of 2019, with  millions of confirmed cases and more than a million deaths \cite{worldometer:...}. On March 11, 2020, COVID-19 was characterized as a pandemic by the World Health Organization (WHO) \cite{WHO}. Due to the huge negative effects of COVID-19, precautionary measures have been aggressively carried out worldwide, such as facial masking, contact tracing, social distancing, and some governmental actions such as lockdowns. Hence, it is significant to analyze the dynamics of COVID-19 so that the effectiveness of those implemented measures can be verified. 

Epidemiological models provide an efficient tool for determining and explaining the dynamics of disease transmission. In the early stage, one of the classical models is the Susceptible-Infectious-Recovery (SIR) model presented by Kermack and McKendrick in 1927 \cite{kermack1927contribution}. This compartmental model computes the theoretical number of susceptible people that became infected by the disease and then recovered. Based on the SIR model, many models have been proposed for different diseases, like SARS and COVID-19 \cite{shi2003stochastic, guanghong2004sars, wang2006spatial}. The Susceptible-Exposed-Infectious-Removed (SEIR) model is proposed to analyze the effect of the precautionary measures on the dynamics of the epidemic \cite{carcione2020simulation}. Taking quarantine into account, a SEIQR model is investigated analytically and numerically \cite{rafiq2020design}. The SEIR model has been extended by two extra classes of populations namely, 'C', which is the number of cumulative cases, and 'D' that is the number of severe, critical, and deceased cases, respectively, to understand the trends of COVID-19 outbreak in Wuhan, China \cite{lin2020conceptual}.  

The SIRD model is governed by three parameters $\beta$, $\gamma$, and $\mu$. Commonly, these parameters are assumed as independent of time during the epidemic. However, the interventions and measures implemented by governments to control the spread of COVID-19 result in significantly varying parameters with time. Thus, if parameters are still considered as constants, the epidemiology model cannot discover the varied dynamic of this disease transmission that is affected by those efficient precautions. Regarding parameters as time-series enable the model to react to the changes in real situations. 

Recently, data-driven algorithms were developed to the determination of such time-varying parameters. Raissi et al. \cite{raissi2019parameter} proposed PINN to approximate parameters, which are constants, of the SIR model using a small data set. Furthermore, Wang et al. \cite{wang2020phase}  divided the period time into four phases so that their parameters become piece-wise constants. Also, there are some studies relating to the time-dependent parameter, like a time-dependent SIR model is proposed to track their parameter time-series \cite{chen2020time}. Similarly, the time-dependent parameters of the SIR model were identified by utilizing a neural network to analyze the COVID-19 spread in South Korea \cite{jo2020analysis}. 


Predicting the trend of COVID-19 transmission is a significant and arduous topic. A model of COVID-19 in China was developed to predict the cumulative number of reported cases \cite{liu2020predicting}. By averaging the slope of approximated parameters of the SIRD model over the last several days, the extrapolated trends of infected, recovered, deceased, and susceptible cases were obtained via solving the model with estimated parameters \cite{magri2020first}. An artificial intelligence method was used to predict the pandemic by training on the 2003 SARS data \cite{yang2020modified}. Recently, several forecasting models including autoregressive integrated moving average, support vector regression, LSTM, and bidirectional LSTM (BiLSTM) were proposed to predict the confirmed cases, deaths, and recoveries \cite{shahid2020predictions}. Similarly, Recurrent Neural Network (RNN), LSTM, BiLSTM, Gated recurrent units (GRUs), and Variational AutoEncoder (VAE) algorithms were employed to forecast the number of infected and recovered cases \cite{zeroual2020deep}.

Based on data collected, the SIRD model is used in this study to simulate the dynamics of the disease. Parameters of this model are regarded as daily and weekly time-varying. As for daily time-varying parameters, it represents that parameters would be varied from each day so that they are more similar to real situations. However, the difficulty of parameter identification increases with the smaller time interval by utilizing some general regression methods such as least square methods. PINN is still able to approximate these time-series parameters with some changes in architecture. Then, they are verified by being substituted into the SIRD model, and then a numerical algorithm is used to solve this ordinary differential equation (ODE) system. 

When it comes to weekly time-series, we divide the collected data week by week and employ the interpolation method for each week to gain more data so that the neural network has enough information to learn. The neural network implemented here is the PINN, which can estimate parameters to fit the collected data \cite{raissi2019physics}. When the identification of these weekly time-series parameters is validated, LSTM is implemented to predict future parameters. In this way, by examining weekly values instead of daily ones, LSTM can still accurately predict future parameters. Another reason to predict the parameters rather than forecast the infectious cases directly is that we could obtain corresponding reproduction number $R_{0}$. It is an essential threshold to reflect the disease transmission speed and verify the effectiveness of these measures in curbing the outbreak of COVID-19. 

\begin{figure}[H]
    \centering
    \includegraphics[scale=0.24]{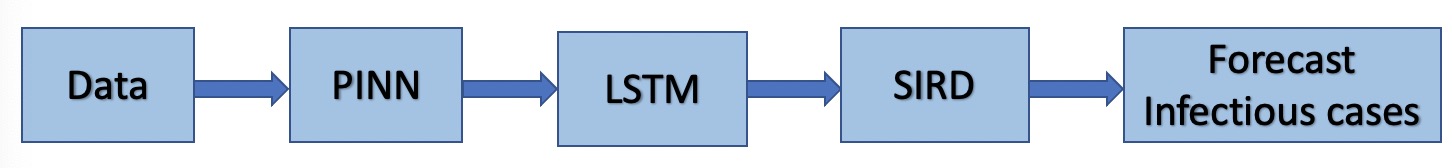}
    \caption{Flow Chart of Data-Driven Algorithm}
    \label{chart}
\end{figure}

In Figure \ref{chart}, we describe that PINN approximates the parameters of the SIRD model by training collected data. Then, LSTM is proposed to predict the future value of parameters. After that, obtained parameters are substituted into the SIRD model so that we could predict the infectious cases. 

\textbf{Outline of the paper:} In Section \ref{epid} we introduce the SIRD model and explain model parameters ($\beta$, $\gamma$, and $\mu$). Following that, we introduce the effective reproduction number and how to compute it. Deep learning neural networks, including feedforward neural network (FNN), PINN, and LSTM, are described briefly in Section 3. In Section 4, the collected data is described, and the loss function of PINN is formulated. After that, we present LSTM neural networks for the future prediction of parameters. The simulation results are discussed in Section 5.

\section{Epidemiology Model}
\label{epid}
In this section, we introduce the SIRD model employed in our study. One promise of the SIRD model is that natural birth and death rates are neglected or equivalent so that the total population is considered as constant. Then, the population could be divided into four mutually exclusive groups, which are susceptible, infected, recovered, and deceased, respectively. Also, the reproduction number is introduced in this section. 

\subsection{SIRD Model} 
We consider the SIRD model described by the following system of ordinary differential equations \cite{ferrari2020modelling}:
\begin{gather}
   \begin{aligned}
    \dot{S}(t) &= -\frac{\beta(t) S(t)I(t)}{N} \\
    \dot{I}(t) &= \frac{\beta(t) S(t)I(t)}{N} - (\gamma(t) + \mu(t)) I(t) \\
    \dot{R}(t) &= \gamma(t) I(t) \\
    \dot{D}(t) &= \mu(t) I(t), 
    \label{sird}
\end{aligned} 
\end{gather}
where $S(t)$, $I(t)$, $R(t)$, and $D(t)$ are the numbers of susceptible, infected, recovered, and deceased individuals, respectively, and $N$ is the total population, $N = S(t) + I(t) + R(t) + D(t)$. Without loss of generality, $N$ is a large number and difficult to compute in some optimization problems. Thus, we used the fraction of the population in each, which means that $S(t)$, $I(t)$, $R(t)$, and $D(t)$ are divided by $N$. In this way, our calculation is simpler, and we can still adhere to the same dynamic system of the epidemic.

The parameter $\beta(t)$ represents the number of contacts each day for infected individuals. Besides, there are two essential assumptions. First, contacts between the infected and uninfected people are sufficient to spread the disease. Second, the population is mixed homogeneously. Thereby, $\beta(t) S(t)$ susceptible people are infected by a patient each day on average, and $\beta(t) S(t)I(t)$ is the number of newly infected people, $\gamma(t)$ is the rate at which infected cases become recovered. Then, there are $\gamma I(t)$ infected individuals turned into the third compartment. Similarly, $\mu(t)$ is the rate of mortality by disease, and $\mu(t) I(t)$ infected individuals are deceased. 

\subsection{Reproduction Number}
There is an essential threshold used to get a better understanding of the transmission rate. This threshold is characterized by the reproduction number, which is an estimate of newly infected cases caused by an infected individual \cite{van2017reproduction}. Assuming there is just one person infected in the beginning, then in the SIRD model, $I_{0} = 1$, $S_{0} = N - 1$, and no recovery at this moment. When we utilized the fraction for each compartment, $I_{0} \approx 0$ and $S_{0} \approx 1$. By the second equation of eq$.$ \eqref{sird}, we have: 
\begin{equation}
    \dot{I}(t) = I(t)[\beta(t) S(t) - (\gamma(t) + \mu(t))] \label{diff1},
\end{equation}
where $\beta(t)$, $\gamma(t)$, and $\mu(t)$ are greater than 0. In eq$.$ \eqref{diff1}, it is simple to obtain $\dot{I} > 0$ when $\beta(t) S(t) > (\gamma(t)+\mu(t))$, which means the number of infected cases keeps increasing. 

A key idea about the reproduction number $R_{0}$ is that it is the gauge of the secondary infections at the beginning of this disease invasion \cite{hatchett2007public}. After the disease invasion, the contact rate $\beta$  decreases because of the reduction of this susceptible group. Those who are infected or have already recovered cannot be infected once again. Thereby, another similar threshold, which is the effective reproduction number $R_{e}$ \cite{ridenhour2018unraveling}, is used to indicate whether the disease keeps spreading, and $R_{0}$ is the upper bound of $R_{e}$ \cite{kergassner2020meso}.
\begin{equation}
    R_{e}(t) = \frac{\beta(t) S(t)}{\gamma(t) + \mu(t)} < \frac{\beta(0)}{\gamma(0) + \mu(0)} = R_{0}.
    \label{re}
\end{equation}

\section{Deep Neural Network}
In this section, we provide a brief description of the deep learning neural network architectures utilized in our algorithm for parameter identification and disease dynamics prediction.

\subsection{Feedforward Neural Network}
The feedforward neural network (FNN) is an artificial neural network with a simple architecture in which the connections between the nodes do not form a cycle. FNN is utilized to approximate the target function by applying several activation functions to input recursively and then output a value or vector that is close to the target values. 

The input layer consists of input neurons, bringing the training data into the network for further processing by hidden layers. Layers between the input layer and output layer are hidden layers, which perform linear or nonlinear transformations of training data to produce an intended output. These transformations are implemented by activation functions. The first hidden layer takes the weighted sums of input through an activation function to produce output for each neuron. Then, these outputs plus a constant, which is bias in the first hidden layer, are viewed as an input of the second hidden layer. After that, as the last layer of the network, the output layer presents desired results \cite{Goodfellow-et-al-2016}.

The forward propagation from one layer to the next layer's nodes is given as follows \cite{gray1992training}.
\begin{equation}
    z_{j}^{l+1} = \sum_{i}^{n_{l}}w^{l}_{ij}f_{l-1}(z_{i}^{l})+b^{l},
    \label{NN1}
\end{equation}

where $n_{l}$ represents the number of neurons in the $l^{th}$ layer, $z_{i}^{l}$, $i = 1,\cdots,n_{l}$, denotes the output of the $i^{th}$ node in $(l-1)^{th}$ layer, $f$ represents the activation function, $w_{ij}^{l}$ are the weights between node $i^{th}$ and node $j^{th}$, $b^{l}$ are the bias in the $l^{th}$ layer. The activation functions employed in this study are the hyperbolic tangent function 

\begin{equation}
    \tanh(x) = \frac{e^{x} - e^{-x}}{e^{x} + e^{-x}} \label{tanh}, 
\end{equation}
and sigmoid function 
\begin{equation}
    \sigma (x) = \frac{1}{1+e^{-x}} \label{sigmoid}.
\end{equation}

When output values are obtained from the output layer, they are used to build a loss function, which is used to estimate the error of the model. An optimizer, like Adam \cite{kingma2014adam} or gradient descent method \cite{amari1993backpropagation}, is employed to minimize the loss function so that the output is close to the observations. According to different problems, the selections of optimizer and loss function are different.   

\subsection{Physics-Informed Neural Network}
The physics-inform neural network (PINN) is a data-driven algorithm to approximate the solution of differential equations and identify parameters \cite{raissi2019physics}. It could utilize any type of neural network architectures, like the FNN and Convolutional neural networks (CNN), as the main framework. The applied activation functions and optimization methods in PINN are the same as the usual deep learning techniques. The fascinating part of this algorithm is the loss function, which is generally comprised of boundary conditions, initial values, and physical constraints. 
 
The outputs of the neural network are constrained to satisfy the system of differential equations by penalizing the residuals of differential equations into the loss function. In these residual equations, the derivatives of outputs with respect to time are computed by a black-box, which is automatic differentiation \cite{lu2019deepxde}. Automatic differentiation is employed to train these models and allows PINN to take derivatives with respect to input coordinates so that it could discover the latent physical laws. More details could be found in \cite{raissi2019physics}.
 
\subsection{Long Short-term Memory Neural Network}
As an advanced tool to process and predict time-series, LSTM is a special type of RNN. The difference between RNN and feedforward neural network is that it is capable of storing past information by taking previous layers' output as an input for the next input. Also, all layers of RNN share the same parameters because they perform the same tasks for the elements of the sequence, which is a so-called recurrent. Based on this feature, RNN could predict sequential information. However, on account of only involving the output of the last layer, it cannot work out long-term memory problems \cite{hochreiter1997long}. Thus, LSTMs are created for solving them. 
\begin{figure}[H]
    \centering
    \includegraphics[scale=0.5]{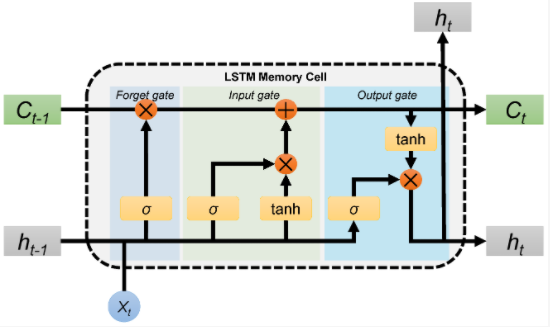}
    \caption{Architecture of LSTM Cell \cite{fan2020comparison}}
    \label{lstm}
\end{figure}
\noindent

LSTMs as shown in Figure \ref{lstm} are chains of memory cells containing three gates, which are input gate, forget gate, and output gate as explained in \cite{fan2020comparison}. The key part of these memory cells is the cell state that could maintain the global information of sequence in each time-step. The content of the cell state would be modified at different time-steps through the forget gates and input gates. More specifically, forget gates decide what information would be thrown away by utilizing the sigmoid function. Instead, input gates decide what new information would be stored by combining sigmoid eq$.$ \eqref{sigmoid} and tanh eq$.$ \eqref{tanh} functions. Then, the cell state is updated and moves forward to the next state by implementing some operations. After that, the output is calculated in the output gate. By analyzing the structure of a memory cell, we could see that LSTMs are capable of storing long-term memory by putting new information into cell states.     

\section{Data and Algorithms}
In this section, we describe the collected data and introduce our algorithms for parameter identification and making predictions. In our work, the parameters of the SIRD model are viewed as daily-varying and weekly-varying parameters to simulate the complex real situation of COVID-19. Based on this perspective, two algorithms are employed to identify daily and weekly time-series parameters. Having the identified parameters by PINN as an input, we introduce the LSTM algorithm for making predictions.

\subsection{Data}
\begin{figure}[H]
    \centering
    \subfloat{\includegraphics[width=2.35in]{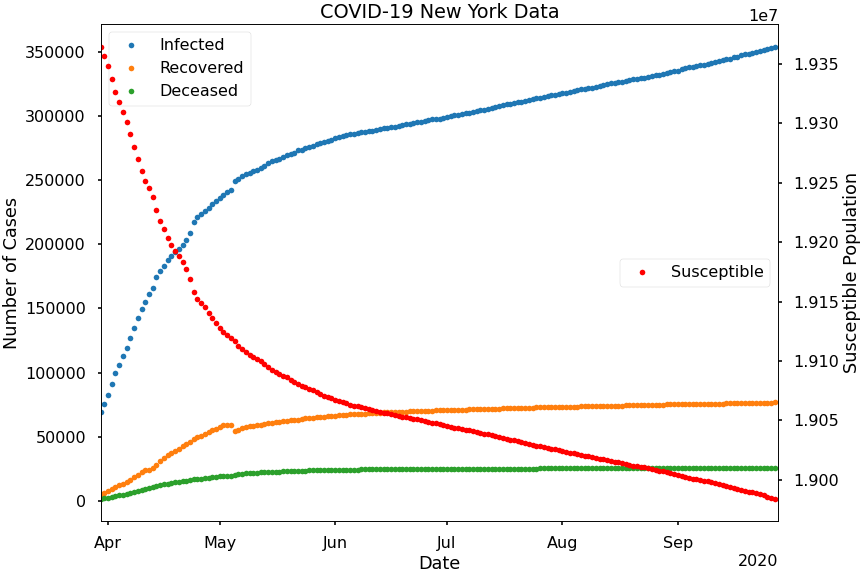}\label{ny_data}}
    \\
    \subfloat{\includegraphics[width=2.35in]{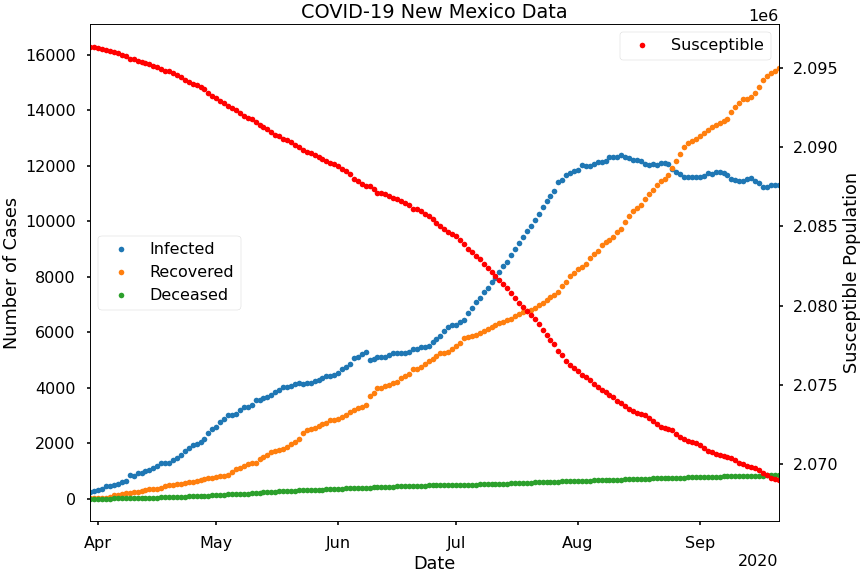}\label{nm_data}}
    \subfloat{\includegraphics[width=2.35in]{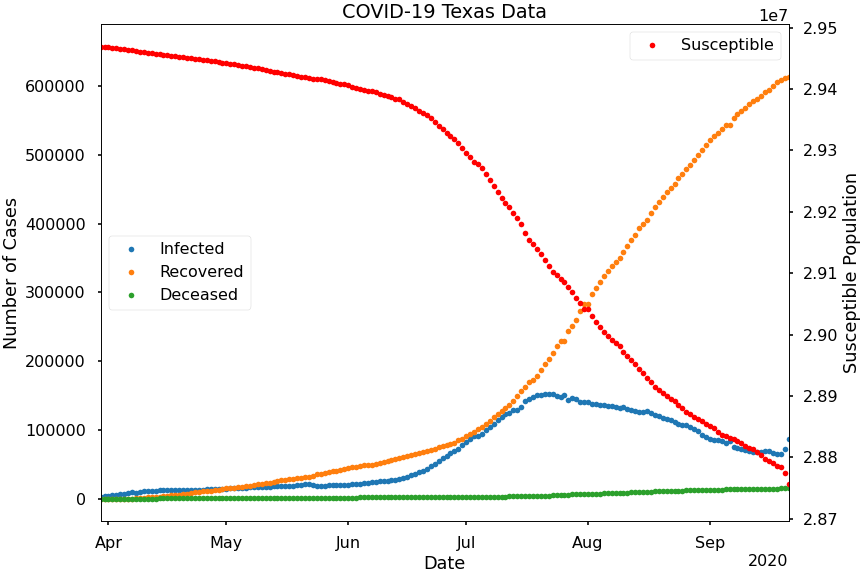}\label{tx_data}}
    \caption{COVID-19 data for New York, New Mexico, and Texas}
    \label{data_source}
\end{figure}

The data considered in our study is downloaded from \url{https://covidtracking.com/data/download} for New York, New Mexico, and Texas states. It is comprised of commutative infected, recovered, and deceased cases for the period from March 30 to September 30. A plot of the data is shown in Figure \ref{data_source}.

As can be observed, the three states have different dynamics between the compartments. In New York, while the number of cases in I, R, and D compartments is increasing, the number of recovered and deceased cases is relatively small compared to infections. New Mexico data shows that the infected and recovered cases kept almost the same increasing rate before August but then the infections started decreasing slowly, while the number of recovered individuals kept increasing at a high rate. As for Texas, the infected population started growing fast from the second half of June, and then went down after about 30 days. 

\subsection{Parameter Identification Algorithm}
In this subsection, we describe how to identify parameters of the SIRD model by utilizing PINN. 

\subsubsection{PINN Architecture for SIRD Model}
The utilized PINN for learning the parameters is composed of an FNN architecture. This architecture consists of 4 hidden layers and 60 neurons in each layer. The activation functions of hidden layers and the output layer are the hyperbolic tangent function, and sigmoid function, respectively. 

 \begin{figure}[H]
    \centering
    \includegraphics[scale=0.22]{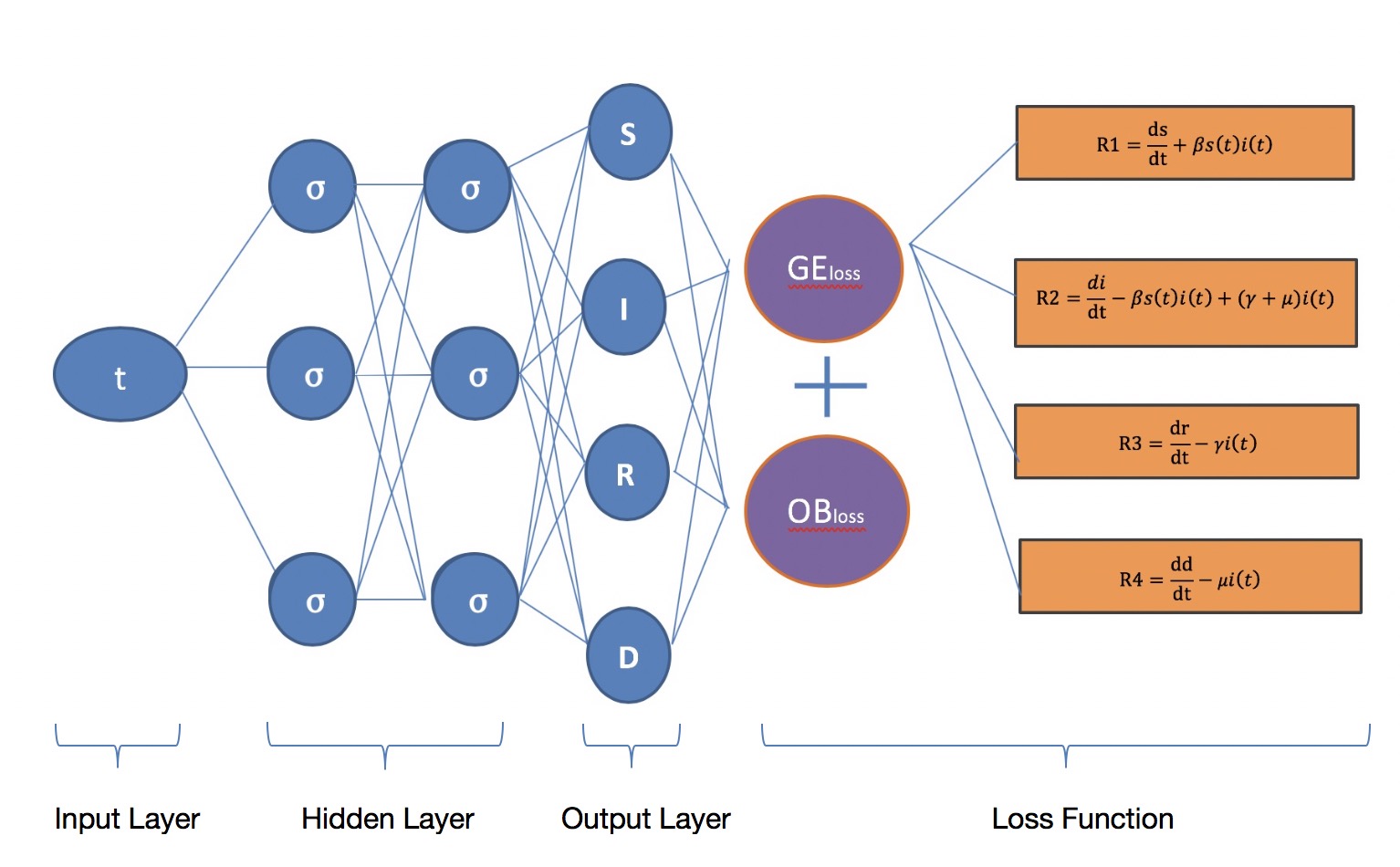}
    \caption{Physics Informed Neural Network}
    \label{pinn_sird}
\end{figure}

Figure \ref{pinn_sird} shows the basic framework, which is FNN, and the loss function of PINN in our study. The loss function of PINN has two parts, namely $GE_{loss}$ and $OB_{loss}$ in this figure. The $GE_{loss}$ represents the residuals that are obtained from the SIRD model eq$.$ \eqref{sird} by subtracting the right side from the left side:
 \begin{gather}
    \begin{aligned}
    R_{1} &= \frac{dS}{dt} + \beta(t) S(t)I(t) \\
    R_{2} &= \frac{dI}{dt} - \beta(t) S(t)I(t) + (\gamma(t)+\mu(t))I(t) \\
    R_{3} &= \frac{dR}{dt} - \gamma(t) I(t) \\
    R_{4} &= \frac{dD}{dt} - \mu(t) I(t).
    \end{aligned}
    \label{error func}
 \end{gather}

The $OB_{loss}$ is the mean squared error between the outputs of the neural network and data. We employ Adam algorithm, which is the first-order gradient-based optimization of stochastic objective functions \cite{kingma2014adam}, to update the parameters of the neural network by minimizing the loss function. 

\subsubsection{Parameters Identification Algorithm}
Algorithm \ref{daily algo} describes how to use PINN to identify the daily time-varying parameters. The input is time $t$, and outputs include the three parameters and the four compartments of the SIRD model. The weights $w$ and bias $b$ are initialized randomly.

\begin{algorithm}[H]
\KwData{$S$,$I$,$R$,$D$,$t$}
\emph{Initialize weights $w$, bias $b$}\;
\For{number of epochs}{
Each compartment of SIRD is obtained from forwarding propagation of neural network
\begin{equation*}
    \hat{S},\hat{I},\hat{R},\hat{D},\beta,\gamma,\mu= NN(t)
\end{equation*}

Build the loss function components, where $N_{ob}$ is the number of collected data of each compartment, and $N_{f}$ is the number of collocation points
\begin{equation*}
OB_{loss} = \frac{1}{N_{ob}}\sum_{i=1}^{N_{ob}}\Big(\big|\hat{S}^{i}-S^{i}\big|+\big|\hat{I}^{i}-I^{i}\big|+\big|\hat{R}^{i}-R^{i}\big| +\big|\hat{D}^{i}-D^{i}\big|\Big)\;
\end{equation*} 
$OB_{loss}$ denotes the difference between the output of the neural network and observation data
\begin{multline*}
 GE_{loss} = \frac{1}{N_{t}}\sum_{i=1}^{N_{t}}\Big(\Big|\frac{d\hat{S}^{i}}{dt^{i}} + \beta \hat{S}^{i}\hat{I}^{i}\Big| + \Big|\frac{d\hat{I}^{i}}{dt^{i}} - \beta \hat{S}^{i}\hat{I}^{i} + (\gamma + \mu\hat{I}^{i})\Big|\\
     + \Big|\frac{d\hat{R}^{i}}{dt^{i}} -\gamma \hat{I}^{i}\Big| + \Big|\frac{d\hat{D}^{i}}{dt^{i}} - \mu\hat{D}^{i}\Big|\Big)   
\end{multline*}
$GE_{loss}$ represents the SIRD model residual loss, then we have the total loss
\begin{equation*}
loss = OB_{loss} + GE_{loss}
\end{equation*}
Using the Adam optimizer to update the weights and bias by minimizing the loss function
}
\caption{PINN identify daily time-varying parameters}
\label{daily algo}
\end{algorithm}

The approach of identifying the weekly-varying parameters is described in algorithm \ref{weekly algo}. We divide the utilized data of algorithm \ref{daily algo} into weekly intervals, and then they are trained to identify corresponding parameters. Cubic spline interpolation is employed to obtain adequate training set for the neural network training.

\begin{algorithm}[H]
\KwData{$S$,$I$,$R$,$D$,$t$}
\emph{Form the weekly data}\;
\For{each week}{
$S,I,R,D,t = CubicSpline(S,I,R,D,t)$ \;
$N = length(t)$ \Comment*[r]{Each dataset has the same length} 
\emph{Initialize weights $w$, bias $b$,$\beta$,$\gamma$,$\mu$}\;
\For{number of epochs}{
Each compartment of SIRD is obtained from forwarding propagation of neural network
\begin{equation*}
    \hat{S},\hat{I},\hat{R},\hat{D} = NN(t)
\end{equation*}

Build the loss function components.
\begin{equation*}
OB_{loss} = \frac{1}{N}\sum_{i=1}^{N}\Big(\big|\hat{S}^{i}-S^{i}\big|^{2}+\big|\hat{I}^{i}-I^{i}\big|^{2}+\big|\hat{R}^{i}-R^{i}\big|^{2}+\big|\hat{D}^{i}-D^{i}\big|^{2}\Big)\;
\end{equation*} 
$OB_{loss}$ denotes the difference between the output of the neural network and observation data
\begin{multline*}
 GE_{loss} = \frac{1}{N}\sum_{i=1}^{N}\Big(\Big|\frac{d\hat{S}^{i}}{dt^{i}} + \beta \hat{S}^{i}\hat{I}^{i}\Big|^{2} + \Big|\frac{d\hat{I}^{i}}{dt^{i}} - \beta \hat{S}^{i}\hat{I}^{i} + (\gamma + \mu\hat{I}^{i})\Big|^{2}\\
     + \Big|\frac{d\hat{R}^{i}}{dt^{i}} -\gamma \hat{I}^{i}|^{2} + \Big|\frac{d\hat{D}^{i}}{dt^{i}} - \mu\hat{D}^{i}\Big|^{2}\Big)   
\end{multline*}
$GE_{loss}$ represents the SIRD model residual loss
\begin{equation*}
loss = OB_{loss} + GE_{loss}
\end{equation*}
Using the Adam optimizer to update the weights and bias by minimizing the loss function
}
\emph{save the value of $\beta$ $\gamma$ and $\mu$} \\
}
\caption{PINN identify weekly time-varying parameters}
\label{weekly algo}
\end{algorithm}

\subsection{LSTM for Prediction of Infectious Cases}
\begin{figure}[H]
    \centering
    \includegraphics[scale=0.25]{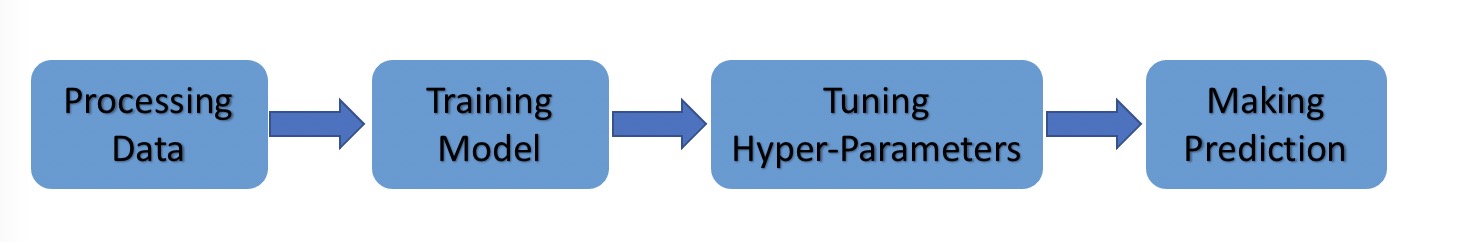}
    \caption{Flow-Chart of LSTM}
    \label{flow-lstm}
\end{figure}

Having the learned parameters by PINN, the LSTM is employed to predict the future parameters by applying Keras as depicted in Figure \ref{flow-lstm}. At first, we normalize these parameters and create new datasets with multiple inputs and one output. Inputs are previous three-time steps parameters, and output, which is the prediction of these inputs, is the next one-time step parameter. Then, the data is used to train the LSTM neural network. We utilize three hidden layers of 80 nodes each to predict the parameters of the next four weeks. 

\section{Simulation}
In this section, we present and discuss our simulation results. The accuracy of the learned parameters by PINN is validated using the relative error in the numerical
solution of System \eqref{sird}, obtained using the Fourth Order Runge-Kutta (RK4) method, with respect to the exact solution. Furthermore, the relative error in the learned solution by PINN with respect to the data is examined. We start by presenting the results for the learned daily and weekly time-varying parameters followed by the associated reproduction numbers. Then, predictions of infectious cases are provided.

\subsection{Time-Varying Parameters}
In general, the daily-varying values of the parameters are more reflective of the daily data than the single constant average values of these parameters. However, identifying the daily-varying parameters is traditionally costly and adds to the complexity of the problem. Fortunately, PINN provides an efficient approach to overcome this difficulty. On the other hand, our simulations show that the learned weekly-varying data produce precisely identified parameters and more stable approximation results. 

Algorithms \ref{daily algo} and \ref{weekly algo} give the identified parameters and the learned values of the SIRD model by PINN. We use the RK4 method to solve the SIRD model with learned parameters. Then, the solutions of the ODE system \eqref{sird} are compared with the learned values.

\begin{figure}[H] 
    \centering
    \subfloat {\includegraphics[width=0.45\textwidth]{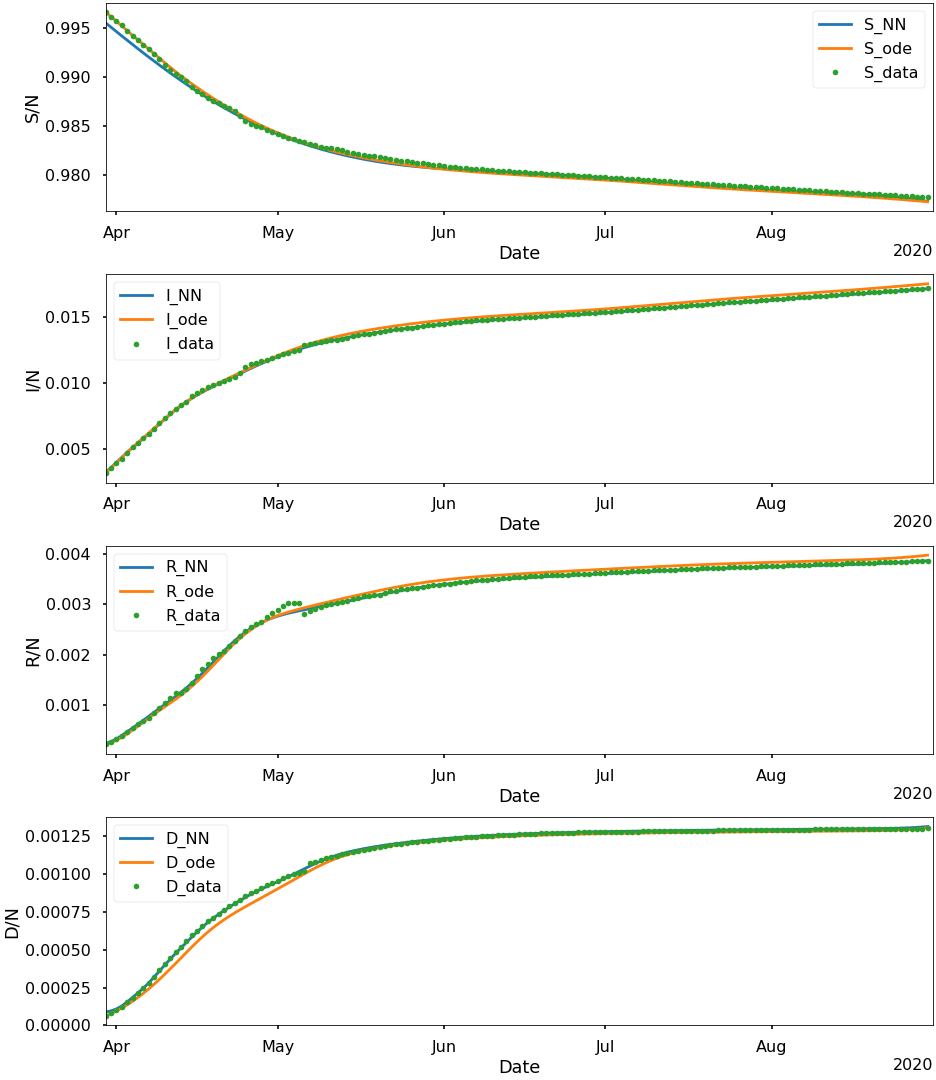}\label{ny_all}}
    \hfill
    \subfloat{\includegraphics[width=0.45\textwidth]{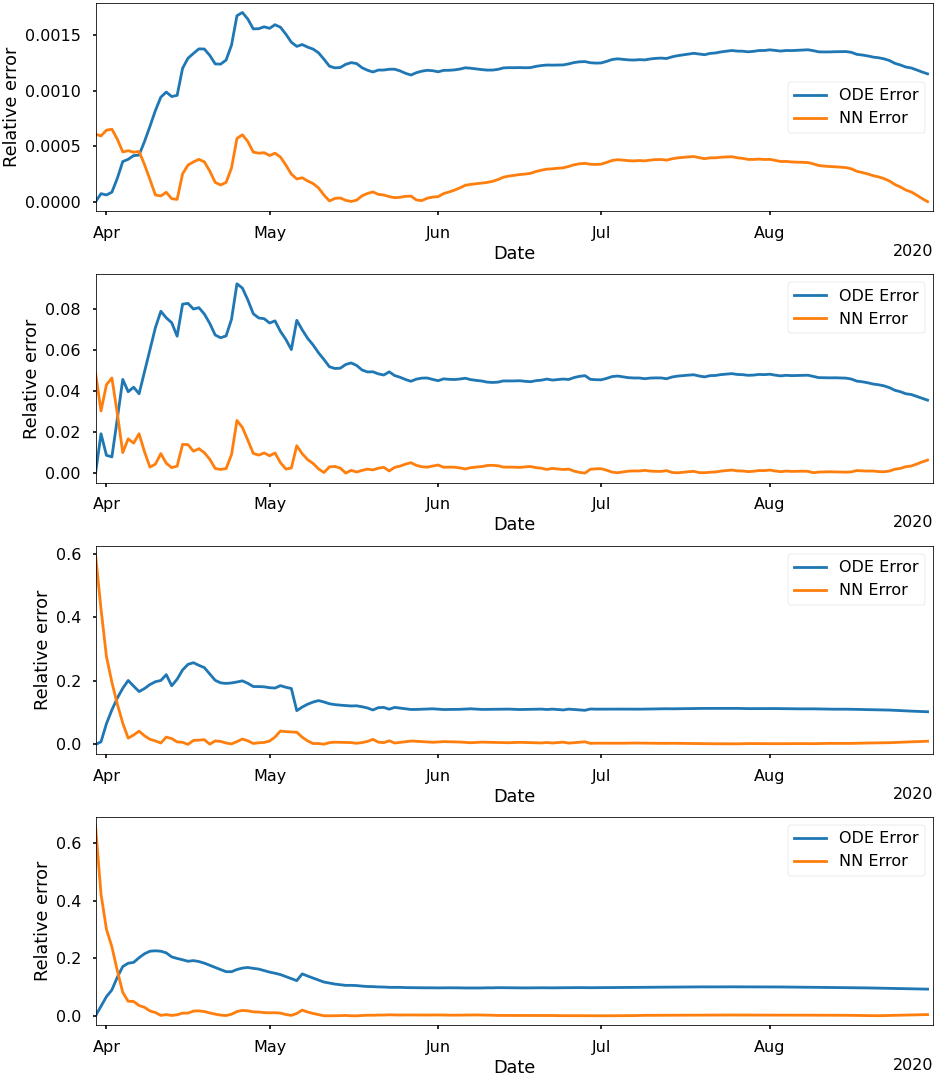}\label{ny_err}}
    \caption{Simulation Results of New York}
    \label{ny_result}
\end{figure}

\begin{figure}[H] 
    \centering
    \subfloat{\includegraphics[width=0.45\textwidth]{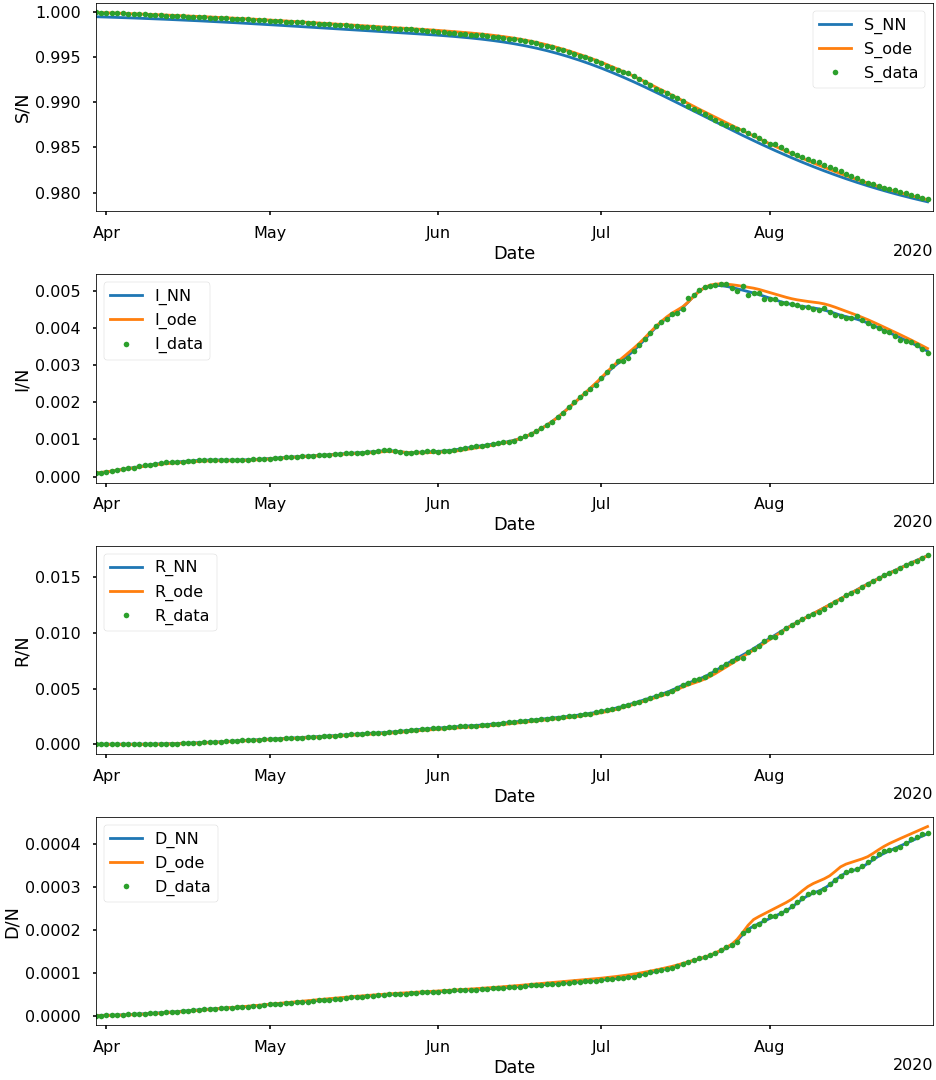}\label{tx_all}}
    \hfill
    \subfloat{\includegraphics[width=0.45\textwidth]{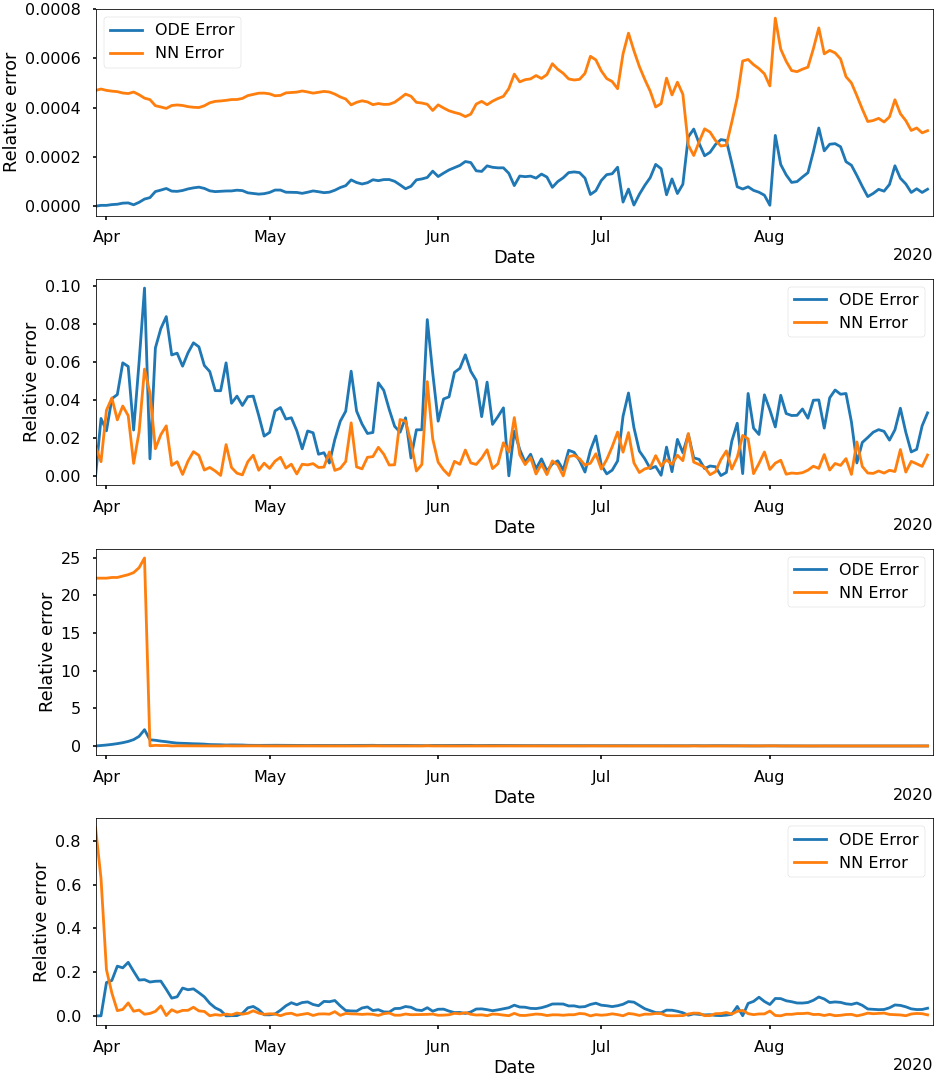}\label{tx_err}}
    \caption{Simulation Results of Texas}
    \label{tx_result}
\end{figure}

\begin{figure}[H] 
    \centering
    \subfloat{\includegraphics[width=0.45\textwidth,keepaspectratio]{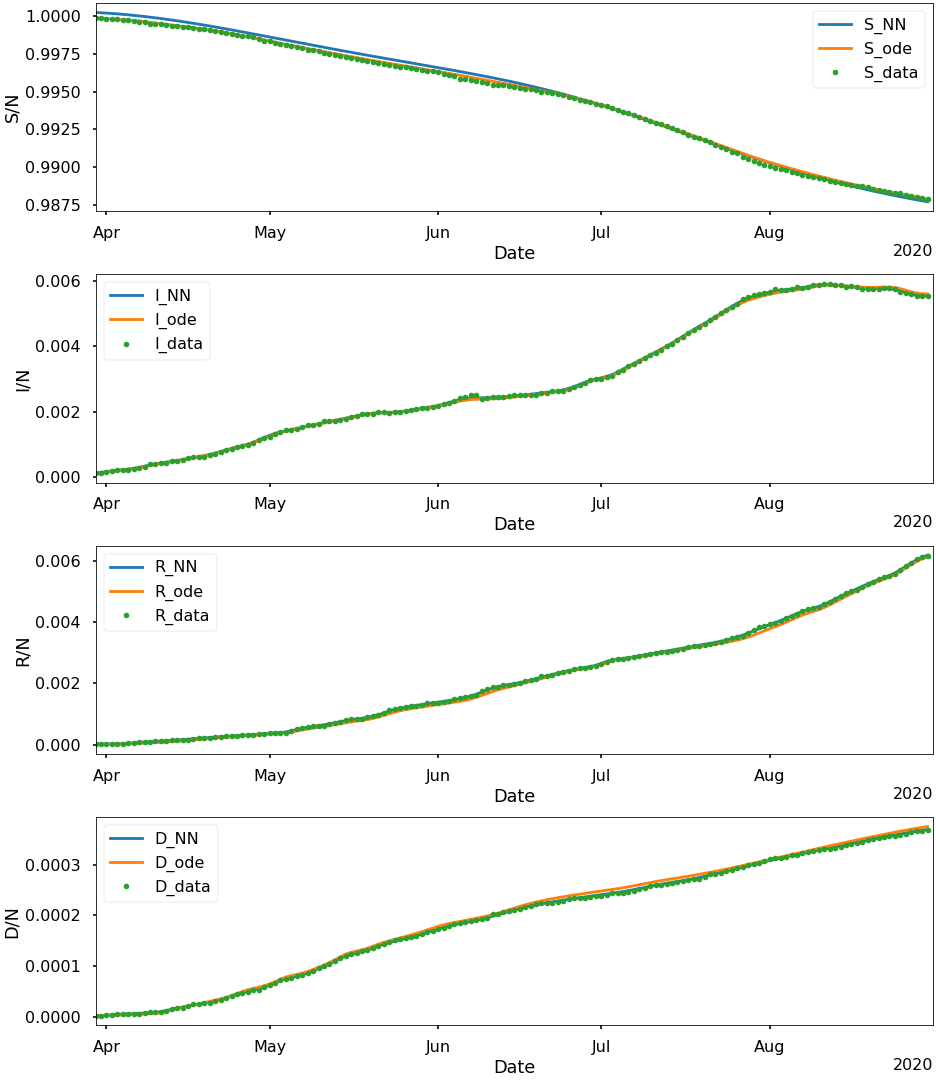}\label{nm_all}}
    \hfill
    \subfloat{\includegraphics[width=0.45\textwidth]{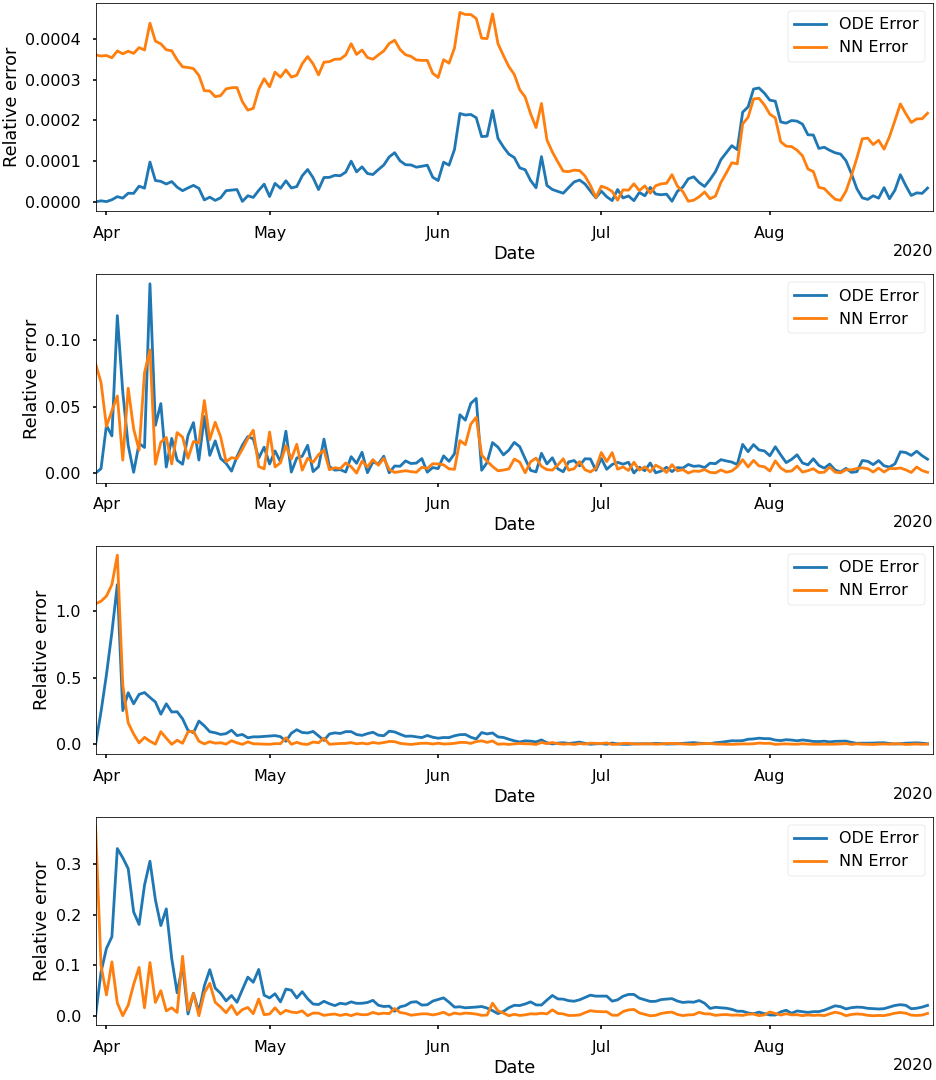}\label{nm_err}}
    \caption{Simulation Results of New Mexico}
    \label{nm_result}
\end{figure}

In Figures \ref{ny_result}, \ref{tx_result}, and \ref{nm_result}, we present the learned values by PINN (NN), the solution of the ODE system (1) (ode), and the reported data (data). Moreover, the relative errors with respect to the reported data are also computed.

\begin{figure}[H] 
    \centering
    \subfloat{\includegraphics[width=0.45\textwidth]{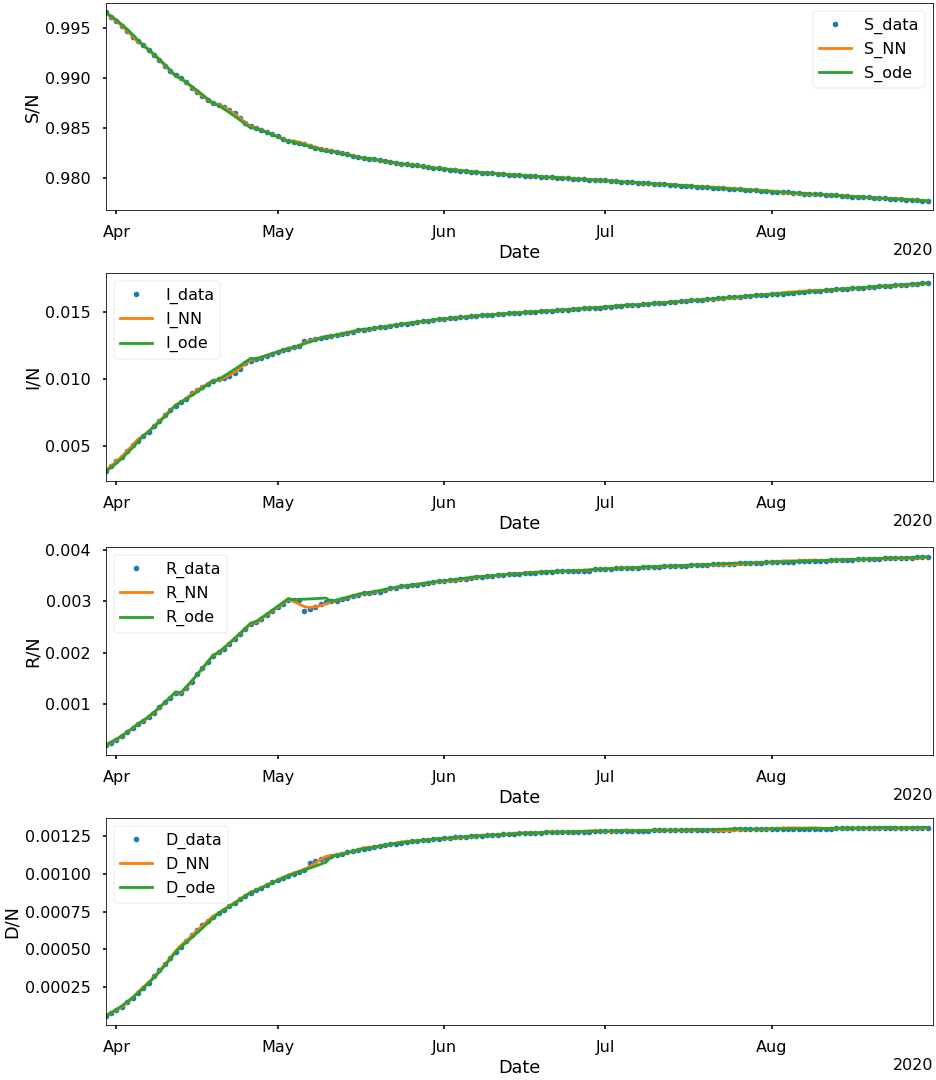}\label{ny_week}}
    \hfill
    \subfloat{\includegraphics[width=0.45\textwidth]{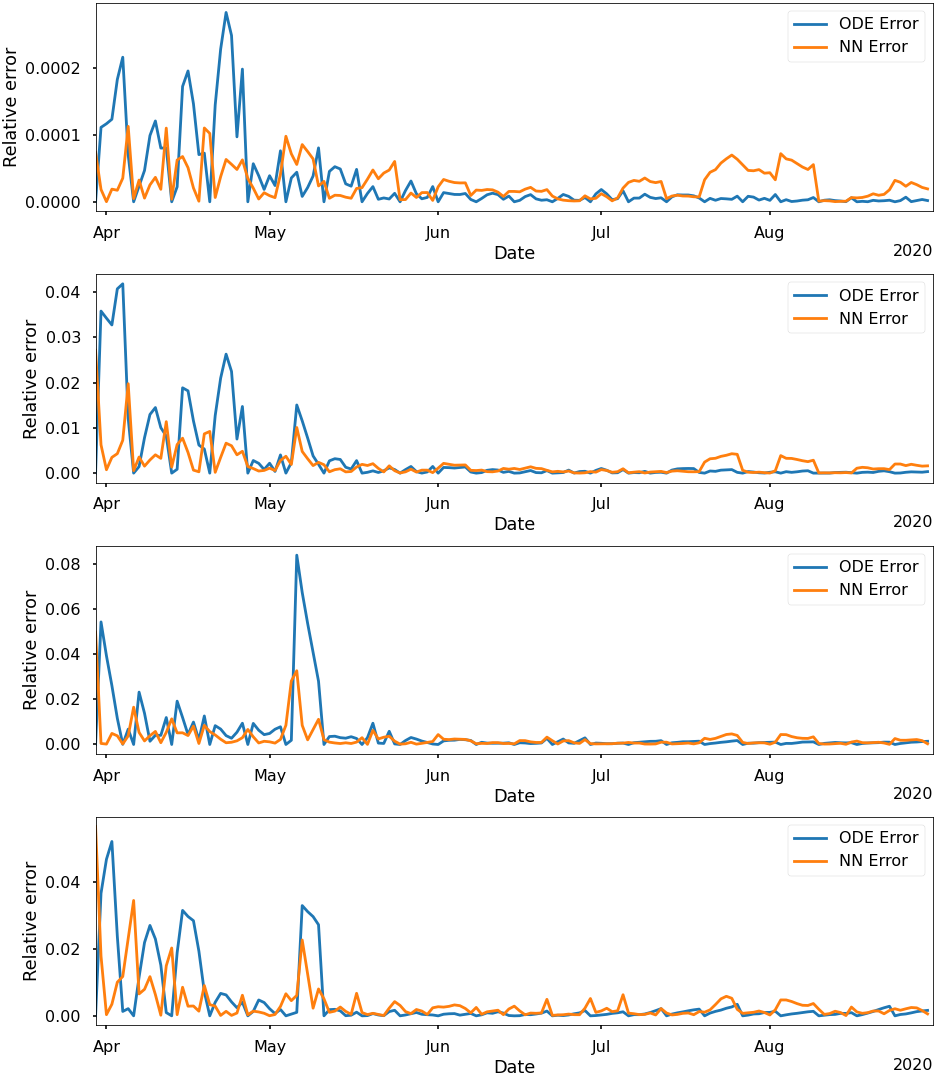}\label{ny_week_err}}
    \caption{Simulation Results of New York}
    \label{ny_week_result}
\end{figure}

\begin{figure}[H] 
    \centering
    \subfloat{\includegraphics[width=0.45\textwidth]{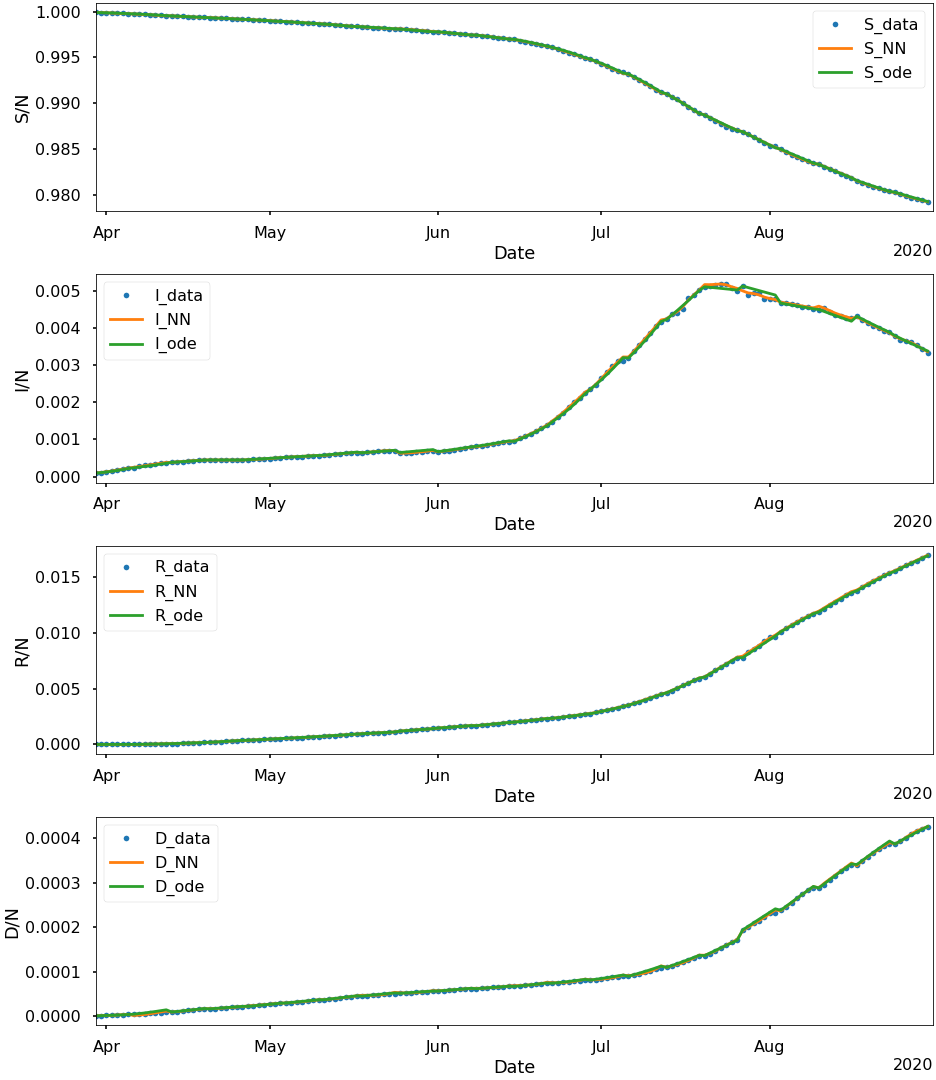}\label{tx_week}}
    \hfill
    \subfloat{\includegraphics[width=0.45\textwidth]{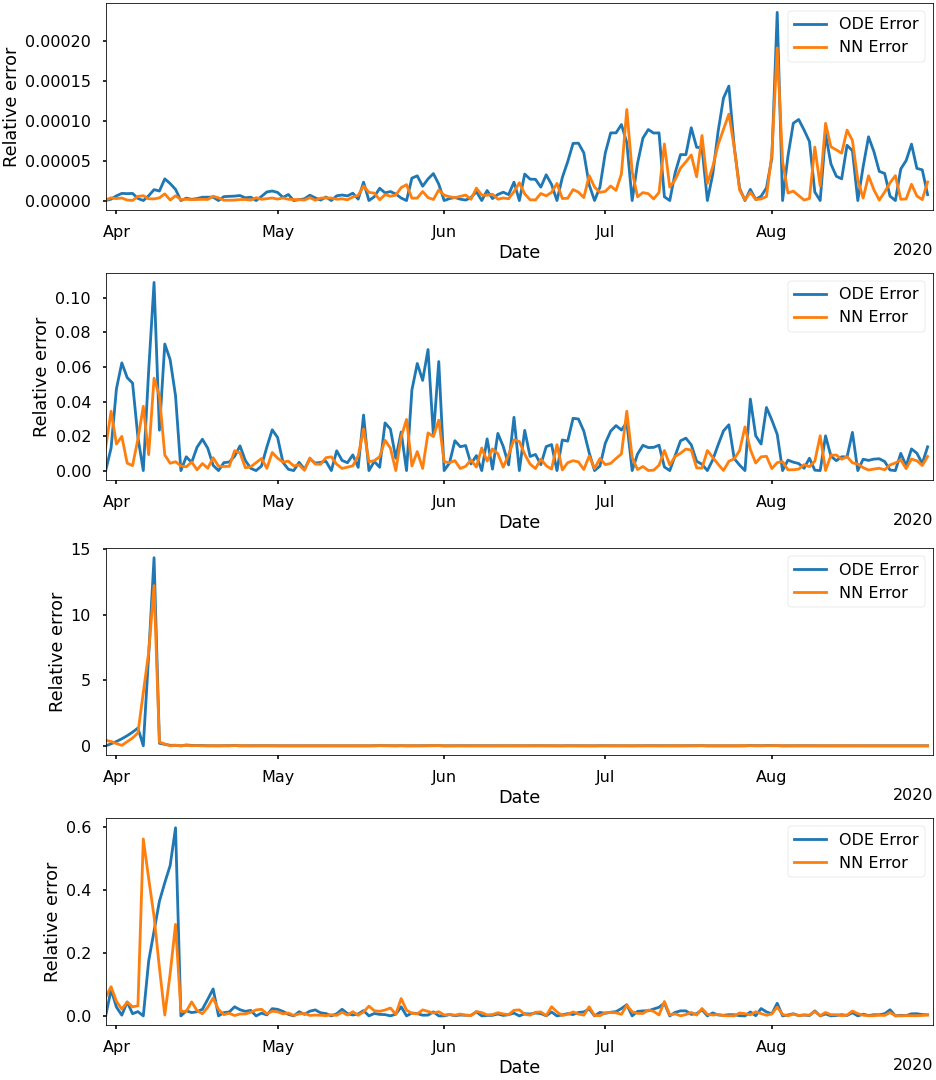}\label{tx_week_err}}
    \caption{Simulation Results of Texas}
    \label{tx_week_result}
\end{figure}

\begin{figure}[H] 
    \centering
    \subfloat{\includegraphics[width=0.45\textwidth]{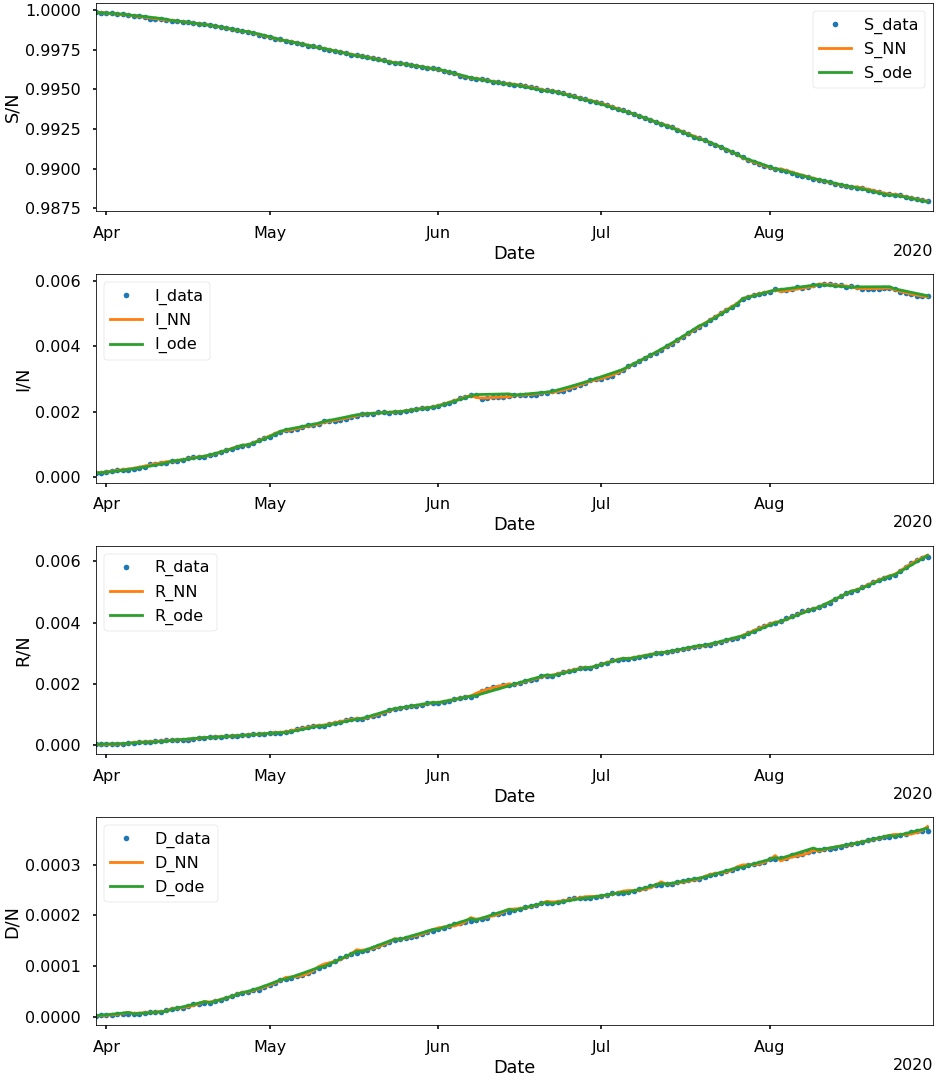}\label{nm_week}}
    \hfill
    \subfloat{\includegraphics[width=0.45\textwidth]{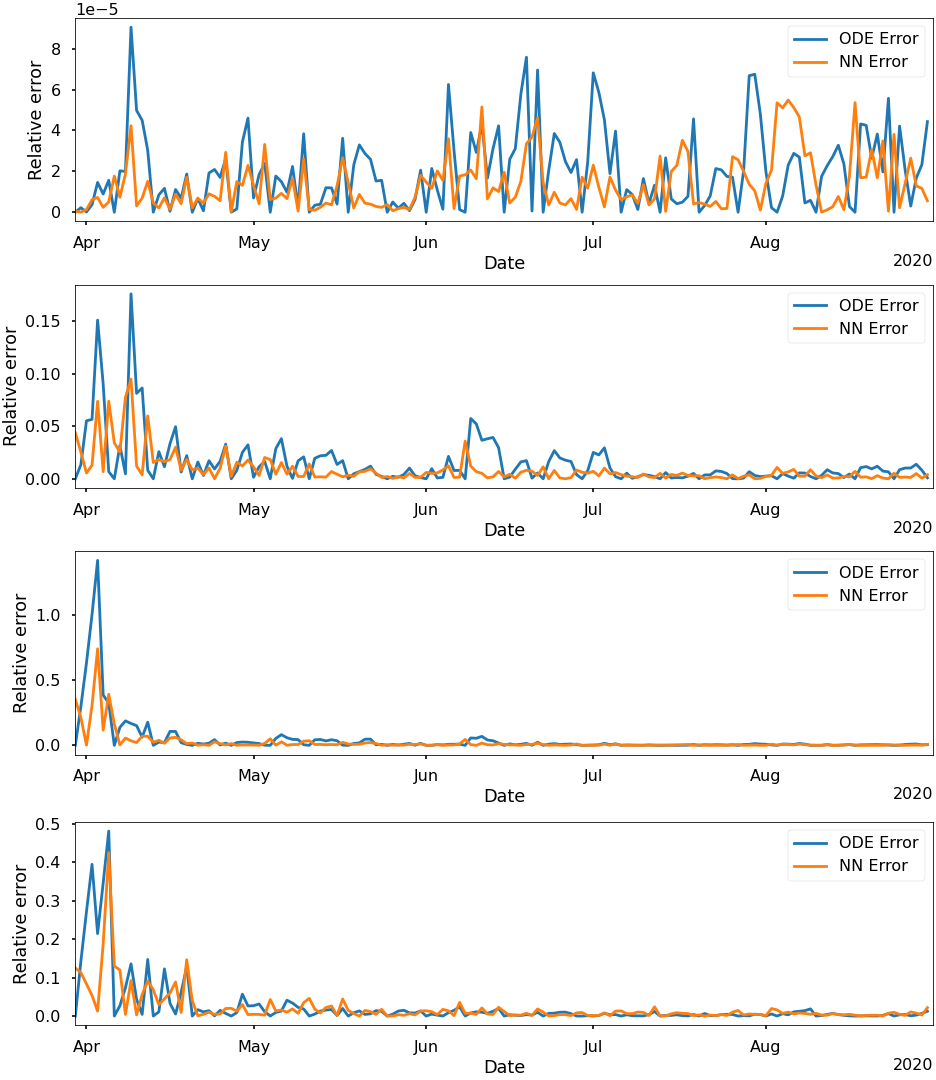}\label{nm_week_err}}
    \caption{Simulation Results of New Mexico}
    \label{nm_week_result}
\end{figure}

There are some large relative errors at the early stage of the outbreak of the disease. This could be due to the small size of the recovered and deceased populations. But the relative errors of the solution of the SIRD model are in a reasonable range, which means that we obtained ideal identified parameters. For the weekly-varying parameters, the period is divided into weeks and the learned results are presented in Figures \ref{ny_week_result}, \ref{tx_week_result}, and \ref{nm_week_result}.


The relative errors in Figures \ref{ny_week_result}, \ref{tx_week_result}, and \ref{nm_week_result} are smaller compared to the results in Figures \ref{ny_result}, \ref{tx_result}, and \ref{nm_result}. Thus, we could infer that the approximation of constant parameters is more accurate and stable. 

\subsection{Simulated Reproduction Numbers}
\begin{figure}[H]
    \centering
    \subfloat[New York]{\includegraphics[width=2.35in]{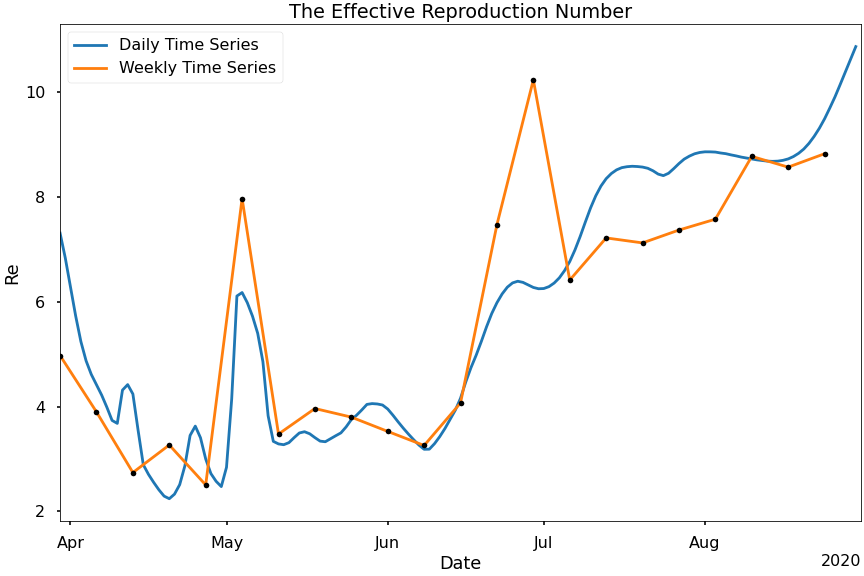}\label{ny_r0}}
    \\
    \subfloat[New Mexico]{\includegraphics[width=2.35in]{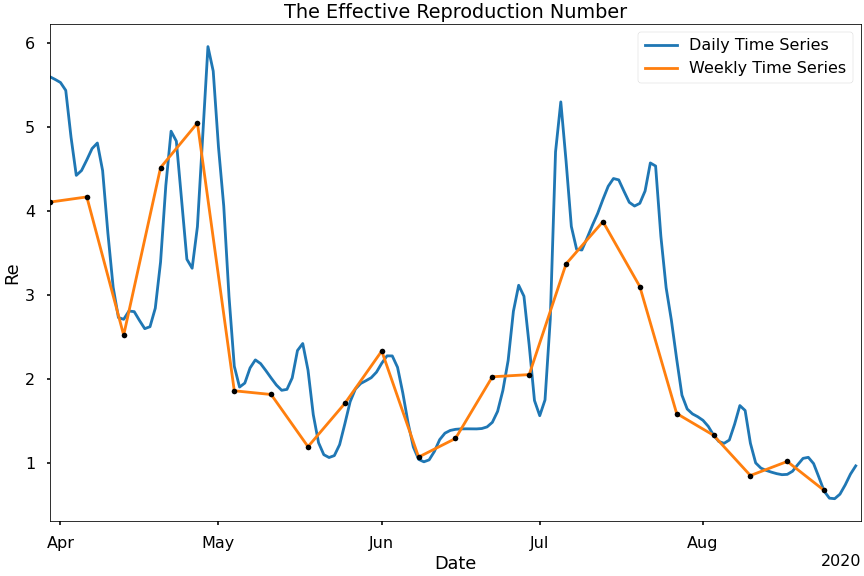}\label{nm_r0}}
    \subfloat[Texas]{\includegraphics[width=2.35in]{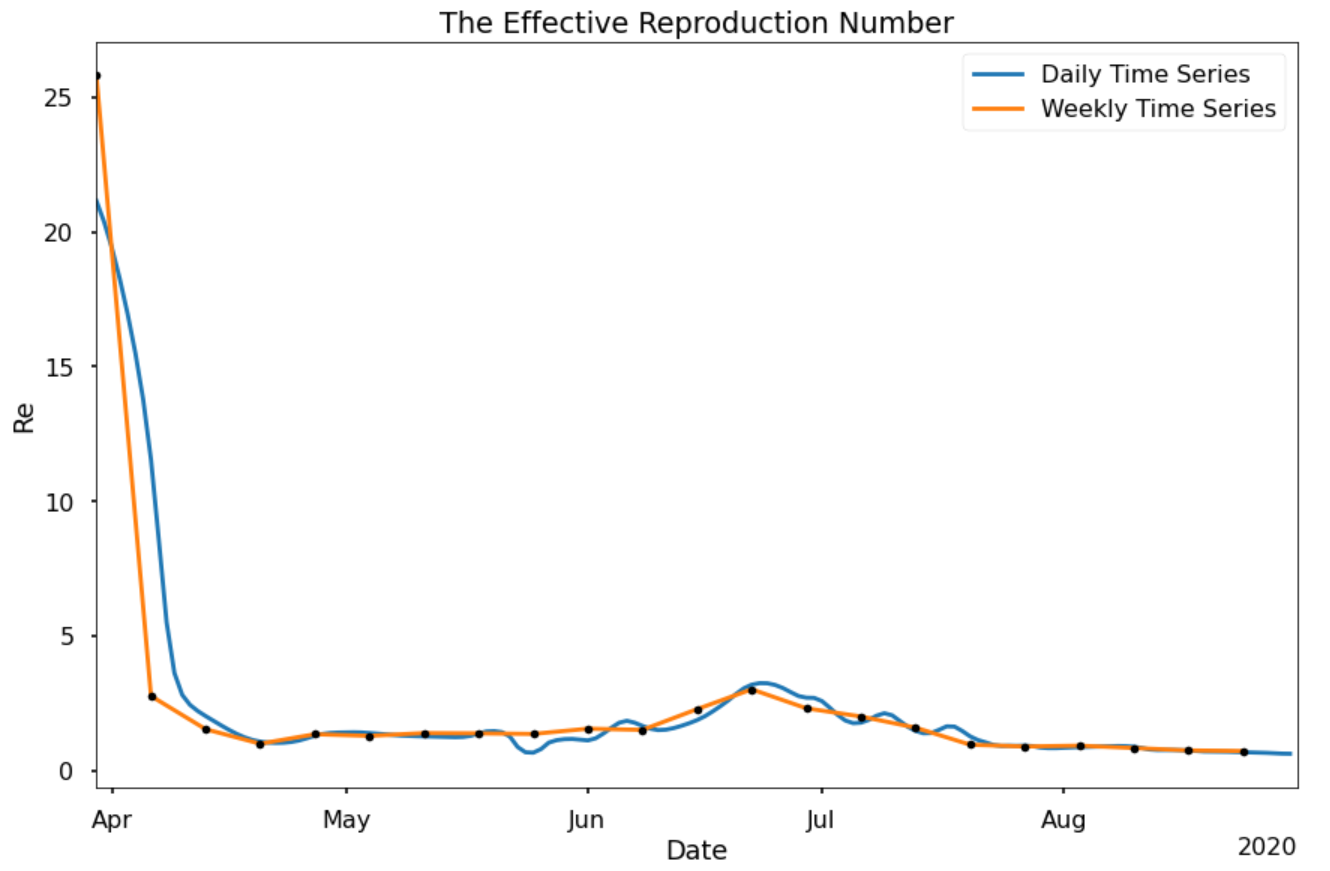}\label{tx_r0}}
    \caption{Effective Reproduction Number}
    \label{r0}
\end{figure}

By eq$.$ \eqref{re}, the effective reproduction number can be obtained using the learned daily and weekly varying parameters. Figure \ref{r0} shows the effective reproduction number in the three states. We could see the daily and weekly time-varying reproduction numbers have the similar trends, which proves the feasibility and effectiveness of our methods to some extent. In the result of New York, we could see the $R_{e}$ is greater than 1 in the whole period, which represents that the number of infectious cases keeps increasing. As for the situation in New Mexico, $R_{e}$ decreases to 1 from August, then the population of infected individuals decreased at that period. By observing the change of $R_{e}$ in Texas, the value is too large in the beginning because the parameter identification in that period is not ideal. But after that period, the $R_{e}$ goes back to normal and oscillates between 0 and 4. Besides, we could see that the growth rate of infected individuals in New York is faster than the other two states by Figure \ref{data_source}. 
\subsection{Prediction of Infectious Cases}
\begin{figure}[H]
    \centering
    \subfloat[New York]{\includegraphics[width=2.35in]{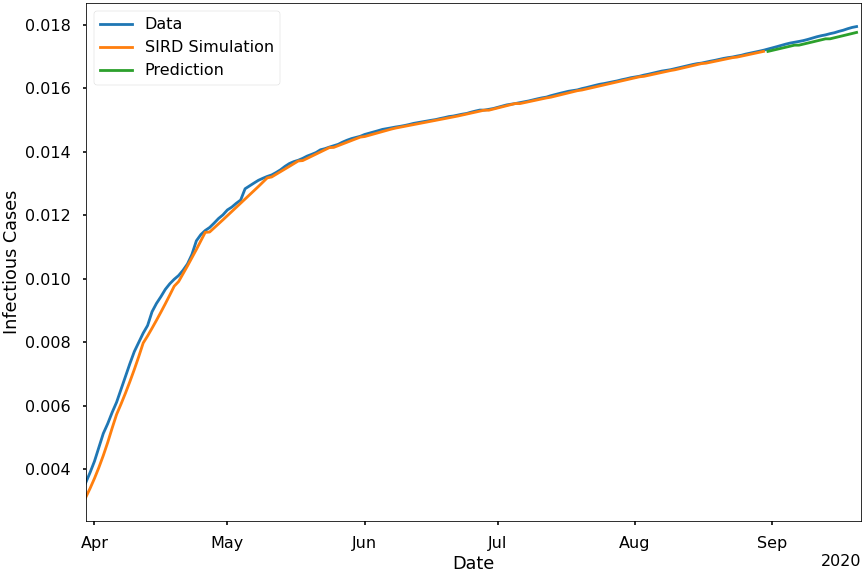}\label{ny_pred}}
    \\
    \subfloat[New Mexico]{\includegraphics[width=2.35in]{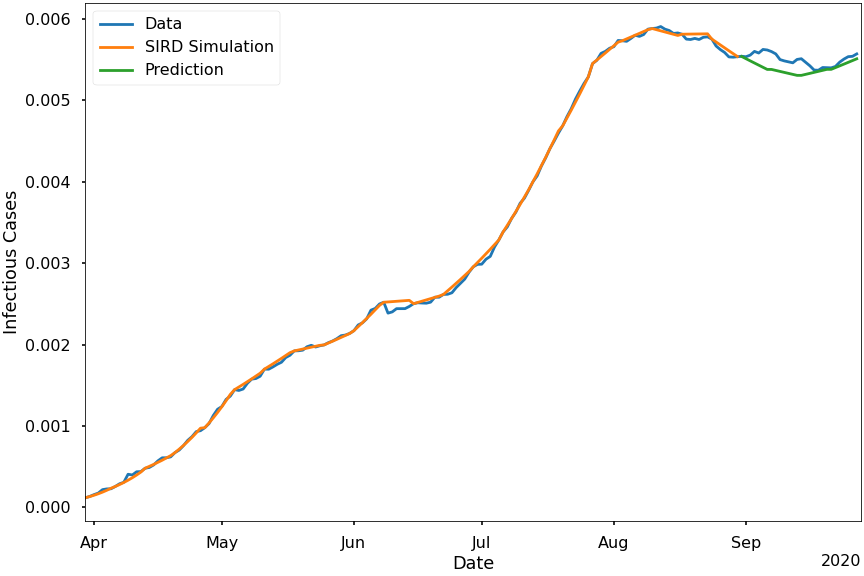}\label{nm_pred}}
    \subfloat[Texas]{\includegraphics[width=2.35in]{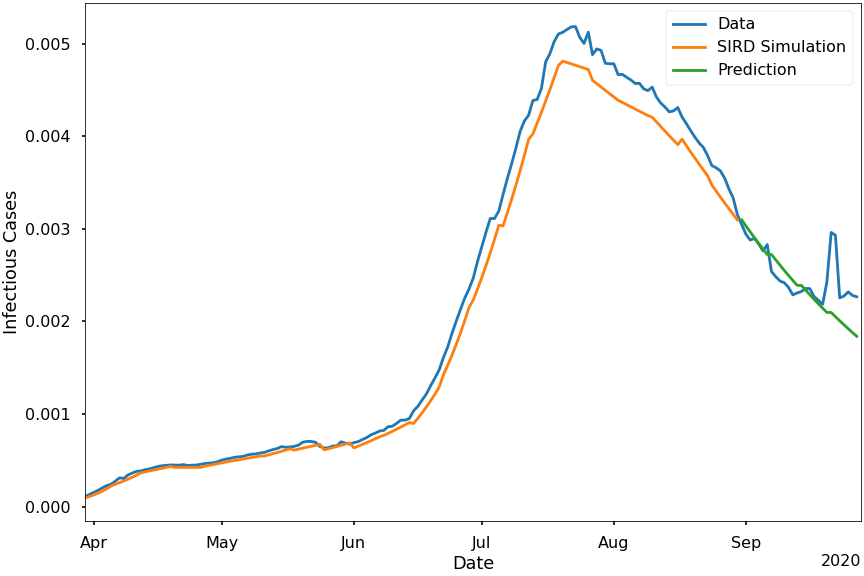}\label{tx_pred}}
    \caption{Predictions}
    \label{prediction}
\end{figure}

The LSTM is employed to predict the future values of the parameters. Figure \ref{prediction} presents the predictions of infected cases for the next four weeks. The prediction of New York is close to the original trend of infectious cases. Predictions in New Mexico show general accord with real situations, but it cannot capture the small fluctuations. As for Texas, there is a sharp increase and decrease in the infectious cases at the end of September, and our prediction could not detect this situation. 

\section{Conclusions}
We introduced a data-driven deep learning approach based on physics informed neural network to solve the epidemiological models and identify weekly and daily time-varying parameters. The PINN was used for parameters identification and solving the ODE system. LSTM was implemented to predict the weekly time-varying parameters and then substituted predicted parameters into the SIRD model so that we could obtain the future trend of COVID-19. The results and errors have shown that the weekly time-varying parameters provided a good fit to the real data. The algorithms could be further developed to achieve a more accurate approximation for complex problems. We intend to explore other architectures of PINN in future work.
\bibliographystyle{tfs}
\bibliography{ref}
\end{document}